\documentclass[11pt]{amsart}

\usepackage{graphicx}
\usepackage[top=26mm, bottom=26mm, left=23mm, right=25mm]{geometry}
\usepackage[dvips]{color}
\usepackage{bm}
\newcommand{\vct}[1]{\bm{#1}}
\newcommand{\mtx}[1]{\mathsf{#1}}

\numberwithin{equation}{section}
\theoremstyle{definition}

\numberwithin{remark}{section}

\numberwithin{definition}{section}

\newcommand{\lsp}{\vspace{3mm}}

\begin{document}

\begin{center}
\textbf{A high-order accurate discretization scheme for variable coefficient elliptic PDEs in the
plane with smooth solutions}

\vspace{2mm}

\textit{P.G. Martinsson, Department of Applied Mathematics, University of Colorado at Boulder}
\end{center}

\vspace{5mm}

\begin{center}
\begin{minipage}{0.8\textwidth}
\noindent\textbf{Abstract:} 
A discretization scheme for variable coefficient elliptic PDEs in the plane is presented.
The scheme is based on high-order Gaussian quadratures and is designed for problems with 
smooth solutions, such as scattering problems involving soft scatterers. The resulting
system of linear equations is very well suited to efficient direct solvers such as
nested dissection and the more recently proposed accelerated nested dissection schemes
with $O(N)$ complexity.
\end{minipage}
\end{center}

\section{Introduction}

\subsection{Background}
This note describes some tentative ideas for how to discretize and solve
a class of variable-coefficient elliptic PDEs with smooth solutions.
The ultimate goal is to device efficient methods for (soft) scattering
problems in the plane, modeled by the Helmholtz' equation
$$
-\Delta\,\phi(x) - \frac{\omega^{2}}{c(x)^{2}}\,\phi(x) = 0,\qquad x \in \mathbb{R}^{2},
$$
where $c(x)$ is a smooth function that is constant outside some domain $\Omega$,
and the ``loading'' of the system is an incoming wave that satisfies the Helmholtz
equation for the constant value of $c$ outside $\Omega$.

The method described is high-order accurate and and leads to a linear system of
algebraic equations that is very well-suited to ``nested-dissection'' type
\textit{direct} (as opposed to \textit{iterative}) solvers.
In a simple implementation, the resulting solver has $O(N^{1.5})$ complexity,
where $N$ is the total number of degrees of freedom in the discretization.
We believe that the cost of the nested dissection step can be further reduced
to $O(N)$ by exploiting techniques similar to those of \cite{2007_leborne_grasedyck,2010_ying_nested_diss,2007_gu_fem}.

In this initial work, we
make several simplifying assumptions in order to investigate the basic viability
of the method. The most significant simplification is that instead of
studying the Helmholtz equation (which has oscillatory solutions), we study the
modified Helmholtz equation (which has non-oscillatory solutions the decay
exponentially fast). The remaining simplifications are, we believe, mostly
cosmetic.

\subsection{Problem statement}
We consider the equation
\begin{equation}
\label{eq:basic}
\left\{\begin{array}{rll}
-\nabla (a(x)\nabla \phi(x)) + b(x)\,\phi(x) =& 0,\qquad &x \in \Omega,\\
                                           \phi_{n}(x) =& v(x),\qquad &x \in \Gamma,
\end{array}\right.
\end{equation}
where $\Omega$ is a box in the plane with boundary $\Gamma = \partial \Omega$,
and where $\phi_{n}$ is the normal derivative of $\phi$.
We assume that the functions $a$ and $b$ are $C^{\infty}$, that $b \geq 0$,
and that for every $x \in \overline{\Omega}$, $a(x)$ is positive definite.
Under these assumptions, the solution $\phi$ and its gradient $\nabla\phi$
will be $C^{\infty}$ in the interior of the domain and can to very high accuracy
be specified via tabulation at Gaussian quadrature nodes. Near the boundary $\Gamma$,
the question of smoothness in general gets complicated, but in this preliminary
report we sidestep this issue by assuming that the given boundary data $v$ is such
that the solution $\phi$ is $C^{\infty}$ on the closed domain $\overline{\Omega}$.
(Ultimately, the technique will be applied to scattering problems and the function
$v$ will be the restriction to $\Gamma$ of the ``incoming wave.'')

\subsection{Outline of the discretization scheme}
\label{sec:sketch_disc}
We propose to tessellate the computational domain into a large number of small squares
and then use the fluxes across the boundaries of each square as the unknown variables in
the model. The flux is represented as a function along the edge; since
this function is smooth, it can very accurately be represented by simply tabulating it
at Gaussian nodes along the edge. We observe that if the fluxes through all four edges
of a box are known, then the values of the potential $\phi$ on all of the edges can
be constructed via the so called ``Neumann-to-Dirichlet'' (N2D) operator for the box. These
operators can cheaply be constructed for all the boxes via a local computation.
Once the N2D operators for all boxes are known, we construct for each edge an equilibrium
equation by combining the N2D operators of the two boxes that share the edge. By combining
the equilibrium equations for all interior edges, we obtain a global equation for all the
interior boundary fluxes. Once this global equation has been solved, the potential on any
box can easily be reconstructed by solving a local Neumann problem on the box (since the
boundary fluxes are now all known).

\subsection{Outline of the linear solver}
\label{sec:sketch_solver}
The discretization described in Section \ref{sec:sketch_disc} results in a large
sparse linear system with a coefficient matrix $\mtx{A}$. If there are $N_{\rm edge}$
edges in the model, and we place $N_{\rm gauss}$ interpolation nodes at each edge, then
$\mtx{A}$ is a block matrix consisting of $N_{\rm edge} \times N_{\rm edge}$ blocks,
each of size $N_{\rm gauss} \times N_{\rm gauss}$.
By ordering the edges in a nested dissection fashion, fill-in can
be limited in the factorization of $\mtx{A}$, resulting in an $O(N_{\rm edge}^{1.5})$ total
cost for the initial solve. (Once one factorization has been executed, subsequent solves
require $O(N_{\rm edge})$ operations.)

The dominant cost in the factorization described above is the inversion or factorization
of dense matrices of size roughly $N_{\rm edge}^{1/2} \times N_{\rm edge}^{1/2}$. These matrices
have internal structure (they are so call \textit{Hierarchically Semi-Separable} (HSS)
matrices) which can be exploited to further reduce the complexity to $O(N_{\rm edge})$.

\section{Discretization}
\label{sec:disc}

We tessellate $\Omega$ into an array of small boxes, and let $\{\Gamma^{(i)}\}_{i \in I_{\rm leaves}}$
denote the collection of \textit{edges} of these boxes, see Figure \ref{fig:domain}.
We include both interior and exterior edges.
For an edge $i$, we define $u^{(i)}$ as the restriction of $\phi$ to $\Gamma^{(i)}$:
$$
u^{(i)}(x) = \phi(x),\qquad\mbox{for }x \in \Gamma^{(i)}.
$$
Further, we define $v^{(i)}$ as the restriction of the normal derivative across $\Gamma^{(i)}$:
$$
v^{(i)}(x) = \left\{\begin{split}
[\partial_{2}\phi](x)\qquad&\mbox{for }x \in \Gamma^{(i)}\mbox{ when } \Gamma^{(i)} \mbox{ is horizontal},\\
[\partial_{1}\phi](x)\qquad&\mbox{for }x \in \Gamma^{(i)}\mbox{ when } \Gamma^{(i)} \mbox{ is vertical},
\end{split}\right.
$$
where we used the short-hand $\partial_{i} = \partial/\partial x_{i}$.

\lsp

\noindent
\fbox{\begin{minipage}{\textwidth}
\textit{\textbf{Observation:}} All the functions $u^{(i)}$ and $v^{(i)}$ are
smooth. They can to very high accuracy be specified by tabulating them at Gaussian
points on the boundary and then interpolate between these points.
\end{minipage}}

\lsp

On each line $\Gamma^{(i)}$, we place $N_{\rm gauss}$ Gaussian nodes.
These points are collected in vectors
$\vct{\gamma}^{(i)} \in \mathbb{R}^{N_{\rm gauss} \times 2}$.
Then we form vectors $\vct{u}^{(i)}, \vct{v}^{(i)} \in \mathbb{R}^{N_{\rm gauss}}$
by collocating the boundary functions $u^{(i)}$ and $v^{(i)}$ at the Gaussian nodes:
\begin{align*}
\vct{u}^{(i)} =&\ u^{(i)}(\vct{\gamma}^{(i)}),\\
\vct{v}^{(i)} =&\ v^{(i)}(\vct{\gamma}^{(i)}).
\end{align*}

\section{The equilibrium equations}
\label{sec:eqm_eqn}

\subsection{Definition of the Neumann-to-Dirichlet operator}
Let $\Omega^{(\tau)}$ be a subdomain of $\Omega$ with edges
$\Gamma^{(i_{1})},\,\Gamma^{(i_{2})},\,\Gamma^{(i_{3})},\,\Gamma^{(i_{4})}$,
as shown in Figure \ref{fig:N2D}. We define the boundary potentials
and boundary fluxes for $\Omega^{(\tau)}$ via
$$
u^{(\tau)} =
\left[\begin{array}{c}
u^{(i_{1})} \\ u^{(i_{2})} \\ u^{(i_{3})} \\ u^{(i_{4})}
\end{array}\right]
\qquad\mbox{and}\qquad
v^{(\tau)} =
\left[\begin{array}{c}
v^{(i_{1})} \\ v^{(i_{2})} \\ v^{(i_{3})} \\ v^{(i_{4})}
\end{array}\right].
$$
Then there exists a unique operator $T^{(\tau)}$
such that with $u^{(\tau)}$ and $v^{(\tau)}$ derived from any solution $\phi$
to (\ref{eq:basic}), we have
\begin{equation}
\label{eq:def_N2D}
u^{(\tau)} = T^{(\tau)}\,v^{(\tau)}.
\end{equation}
This claim follows from the fact that a Neumann boundary value problem on $\Omega^{(\tau)}$
has a unique solution.
The operator $T^{(\tau)}$ is mathematically an integral operator called the \textit{Neumann-to-Dirichlet}
operator.

\lsp

\noindent
\fbox{\begin{minipage}{\textwidth}
\textit{\textbf{Observation:}}
The N2D operator is in general a complicated object. It has a singular kernel
even for domains with smooth boundaries. When the domain boundary has corners,
further complications arise.
However, in our case, all such subtleties can be ignored
since the boundary potentials of interest are all restrictions of functions
that globally solve (\ref{eq:basic}), and are in consequence smooth.
\end{minipage}}

\lsp

The discrete analog  of the equation (\ref{eq:def_N2D}) is
\begin{equation}
\label{eq:def_N2D_disc}
\vct{u}^{(\tau)} = \mtx{T}^{(\tau)}\,\vct{v}^{(\tau)}.
\end{equation}
For our purposes, it is sufficient for the matrix $\mtx{T}^{(\tau)}$ to
correctly construct $\vct{u}^{(\tau)}$ for any \textit{permissible}
vectors $\vct{v}^{(\tau)}$. (By permissible, we mean that they are
the restriction of a function in the solution set under consideration.)
Written out in components, $\mtx{T}^{(\tau)}$ is a $4\times 4$ block matrix that satisfies
\begin{equation}
\label{eq:writtenout}
\left[\begin{array}{c}
\vct{u}^{(i_{1})} \\ \vct{u}^{(i_{2})} \\ \vct{u}^{(i_{3})} \\ \vct{u}^{(i_{4})}
\end{array}\right] =
\left[\begin{array}{cccc}
\mtx{T}^{(\tau,11)} & \mtx{T}^{(\tau,12)} & \mtx{T}^{(\tau,13)} & \mtx{T}^{(\tau,14)} \\
\mtx{T}^{(\tau,21)} & \mtx{T}^{(\tau,22)} & \mtx{T}^{(\tau,23)} & \mtx{T}^{(\tau,24)} \\
\mtx{T}^{(\tau,31)} & \mtx{T}^{(\tau,32)} & \mtx{T}^{(\tau,33)} & \mtx{T}^{(\tau,34)} \\
\mtx{T}^{(\tau,41)} & \mtx{T}^{(\tau,42)} & \mtx{T}^{(\tau,43)} & \mtx{T}^{(\tau,44)}
\end{array}\right]
\left[\begin{array}{c}
\vct{v}^{(i_{1})} \\ \vct{v}^{(i_{2})} \\ \vct{v}^{(i_{3})} \\ \vct{v}^{(i_{4})}
\end{array}\right].
\end{equation}

\subsection{Construction of the N2D operator for a small box}
\label{sec:N2D_leaf}
Let $\Omega^{(\tau)}$ be a small box with edges
$\Gamma^{(i_{1})},\,\Gamma^{(i_{2})},\,\Gamma^{(i_{3})},\,\Gamma^{(i_{4})}$, as shown in Figure \ref{fig:N2D},
and consider the task of constructing a matrix
$\mtx{T}^{(\tau)}$ such that (\ref{eq:def_N2D_disc}) holds for all permissible potentials.
We generate (by brute force) a collection of solutions $\{\phi_{j}\}_{j=1}^{N_{\rm samp}}$
that locally span the solution space to the desired precision. For each $\phi_{j}$, we construct the corresponding
vectors of boundary values
$$
\vct{u}_{j} = \left[\begin{array}{c}
                       \phi_{j}(\vct{\gamma}^{(i_{1})})\\
                       \phi_{j}(\vct{\gamma}^{(i_{2})})\\
                       \phi_{j}(\vct{\gamma}^{(i_{3})})\\
                       \phi_{j}(\vct{\gamma}^{(i_{4})})
                \end{array}\right],
\qquad\mbox{and}\qquad
\vct{v}_{j} = \left[\begin{array}{c}
                        \partial_{2}\phi_{j}(\vct{\gamma}^{(i_{1})})\\
                        \partial_{1}\phi_{j}(\vct{\gamma}^{(i_{2})})\\
                        \partial_{2}\phi_{j}(\vct{\gamma}^{(i_{3})})\\
                        \partial_{1}\phi_{j}(\vct{\gamma}^{(i_{4})})
                        \end{array}\right]
$$
and then we construct via a least squares procedure a matrix $\mtx{T}^{(\tau)}$ such that the equation
\begin{equation}
[\vct{u}_{1}\ \vct{u}_{2}\ \cdots\ \vct{u}_{N_{\rm samp}}] =
\mtx{T}^{(\tau)}\,
[\vct{v}_{1}\ \vct{v}_{2}\ \cdots\ \vct{v}_{N_{\rm samp}}]
\end{equation}
holds to within the specified tolerance $\varepsilon$.

The sample functions $\phi_{j}$ are constructed by solving a set of local problems
on a patch $\Psi$ that covers the domain $\Omega^{(\tau)}$, as shown in Figure \ref{fig:local_patches}.
The local problems read
\begin{equation}
\label{eq:patch_equation}
\left\{\begin{array}{rll}
-\nabla (\tilde{a}(x)\nabla \phi_{j}(x)) + \tilde{b}(x)\,\phi_{j}(x) =& 0,\qquad &x \in \Psi,\\
                                           \partial_{n}\phi(x) =& v_{j}(x),\qquad &x \in \partial \Psi,
\end{array}\right.
\end{equation}
where $\tilde{a}$ and $\tilde{b}$ are functions chosen so that:
\begin{enumerate}
\item For $x \in \Omega^{(\tau)}$, we have $\tilde{a}(x) = a(x)$ and $\tilde{b}(x) = b(x)$.
\item The equation (\ref{eq:patch_equation}) is easy to solve.
\end{enumerate}
The collection of boundary data $\{v_{j}\}_{j=1}^{N_{\rm samp}}$ is chosen so that the solution
space is sufficiently ``rich.''

\subsection{Assembling a global equilibrium equation}
In this section, we will formulate a linear equation that relates the following variables:
\begin{align*}
\mbox{Given data:}&\hspace{5mm} \{\vct{v}^{(i)}\,\colon\,i\mbox{ is an edge that is exterior to }\Omega\},\\
\mbox{Sought data:}&\hspace{5mm} \{\vct{v}^{(i)}\,\colon\,i\mbox{ is an edge that is interior to }\Omega\}.
\end{align*}
Let $N_{\rm edge}$ denote the number of interior edges. Then the coefficient matrix of the linear system
will consist of $N_{\rm edge} \times N_{\rm edge}$ blocks, each of size $N_{\rm gauss} \times N_{\rm gauss}$.
Each block row in the system will have at most $7$ non-zero blocks. To form this matrix, let $i$ denote
an interior edge. Suppose for a moment that $i$ is a vertical edge. Let $\tau_{1}$ and $\tau_{2}$ denote
the two boxes that share the edge $i$, let $\{m_{1},\,m_{2},\,m_{3},\,m_{4}\}$ denote the edges of $\tau_{1}$,
and let $\{n_{1},\,n_{2},\,n_{3},\,n_{4}\}$ denote the edges of $\tau_{2}$, see Figure \ref{fig:merge_new}.
The N2D operator for $\tau_{1}$ provides an equation for the boundary fluxes of the left box:
\begin{equation}
\label{eq:eqn_for_tau1}
\vct{u}^{(m_{2})} =
\mtx{T}^{(\tau_{1},21)}\,\vct{v}^{(m_{1})} +
\mtx{T}^{(\tau_{1},22)}\,\vct{v}^{(m_{2})} +
\mtx{T}^{(\tau_{1},23)}\,\vct{v}^{(m_{3})} +
\mtx{T}^{(\tau_{1},24)}\,\vct{v}^{(m_{4})}.
\end{equation}
Analogously, the N2D operator for $\tau_{2}$ provides the equation
\begin{equation}
\label{eq:eqn_for_tau2}
\vct{u}^{(n_{4})} =
\mtx{T}^{(\tau_{2},41)}\,\vct{v}^{(n_{1})} +
\mtx{T}^{(\tau_{2},42)}\,\vct{v}^{(n_{2})} +
\mtx{T}^{(\tau_{2},43)}\,\vct{v}^{(n_{3})} +
\mtx{T}^{(\tau_{2},44)}\,\vct{v}^{(n_{4})}.
\end{equation}
Observing that $m_{2} = n_{2} = i$, we see that $\vct{u}^{(m_{2})} = \vct{u}^{(n_{4})}$, and consequently
(\ref{eq:eqn_for_tau1}) and (\ref{eq:eqn_for_tau2}) can be combined to form the equation
\begin{multline}
\label{eq:eqn_for_gamma_i}
\mtx{T}^{(\tau_{1},21)}\,\vct{v}^{(m_{1})} +
\mtx{T}^{(\tau_{1},22)}\,\vct{v}^{(i)} +
\mtx{T}^{(\tau_{1},23)}\,\vct{v}^{(m_{3})} +
\mtx{T}^{(\tau_{1},24)}\,\vct{v}^{(m_{4})}\\ =
\mtx{T}^{(\tau_{2},41)}\,\vct{v}^{(n_{1})} +
\mtx{T}^{(\tau_{2},42)}\,\vct{v}^{(n_{2})} +
\mtx{T}^{(\tau_{2},43)}\,\vct{v}^{(n_{3})} +
\mtx{T}^{(\tau_{2},44)}\,\vct{v}^{(i)}.
\end{multline}
The collection of all equations of the form (\ref{eq:eqn_for_gamma_i}) for interior vertical edges,
along with the analogous set of equations for all interior horizontal edges forms the global equilibrium
equation.

\section{Efficient direct solvers}

This section describes a direct solver for the global equilibrium
equation constructed in Section \ref{sec:eqm_eqn}. The idea is to
partition the box $\Omega$ into a quad-tree of boxes, and then to
construct the N2D operator $\mtx{T}^{(\tau)}$ for each box $\tau$
in the tree. The first step is to loop over all leaf nodes of the
tree and construct the N2D operator via the procedure described
in Section \ref{sec:N2D_leaf}. Then we execute an upwards sweep
through the tree, where we construct the N2D operator for a box
by merging the operators of its four children. For simplicity
(and also computational efficiency) we execute each ``merge-four''
operation as a set of three ``merge-two'' operations.

Let $N$ denote the size of the coefficient matrix. Then Section \ref{sec:slow_merge}
describes a procedure with $O(N^{1.5})$ complexity, and Section \ref{sec:fast_merge}
sketches out how the procedure can be accelerated to $O(N)$ complexity. Before
describing the fast solvers, we describe a hierarchical decomposition of the domain
in Section \ref{sec:quadtree}.

\subsection{A quad-tree on the domain}
\label{sec:quadtree}
A standard quad-tree is formed on the computational domain $\Omega$
as follows: Let $\Omega^{(1)} = \Omega$ be the \textit{root} of the tree,
as shown in Figure \ref{fig:tree}(a).
Then split $\Omega^{(1)}$ into four disjoint boxes
$$
\Omega^{(1)} = \Omega^{(2)} \cup \Omega^{(3)} \cup \Omega^{(4)} \cup \Omega^{(5)},
$$
as shown in Figure \ref{fig:tree}(b). Continue by splitting each of the four
boxes into four smaller equisized boxes:
\begin{align*}
\Omega^{(5)} =&\ \Omega^{( 6)} \cup \Omega^{( 7)} \cup \Omega^{( 8)} \cup \Omega^{( 9)},\\
\Omega^{(6)} =&\ \Omega^{(10)} \cup \Omega^{(11)} \cup \Omega^{(12)} \cup \Omega^{(13)},\\
\Omega^{(7)} =&\ \Omega^{(14)} \cup \Omega^{(15)} \cup \Omega^{(16)} \cup \Omega^{(17)},\\
\Omega^{(8)} =&\ \Omega^{(18)} \cup \Omega^{(19)} \cup \Omega^{(20)} \cup \Omega^{(21)},
\end{align*}
as shown in Figure \ref{fig:tree}(c). The process continues until each box is small
enough that the N2D operator for each leaf can easily be constructed via the procedure
described in Section \ref{sec:N2D_leaf}. The levels of the tree are ordered so that $\ell = 0$
is the coarsest level (consisting only of the root), $\ell = 1$ is the level with four boxes, etc.
We let $L$ denote the total number of levels in the tree.

\subsection{Simple construction of the N2D operator for a parent}
\label{sec:slow_merge}

Suppose that $\sigma$ is a box with children $\nu_{1}$ and $\nu_{3}$ as shown in Figure \ref{fig:merge},
and that we know the matrices $\mtx{T}^{(\nu_{1})}$ and $\mtx{T}^{(\nu_{3})}$ associated with the
children. We seek to construct the matrix $\mtx{T}^{(\sigma)}$. The equilibrium equations
for the two children read
\begin{eqnarray}
\label{eq:ida1a}
\vct{u}^{(m_{i})} =& \sum_{j=1}^{4} \mtx{T}^{(\nu_{1},ij)}\,\vct{v}^{(m_{j})},\qquad i = 1,\,2,\,3,\,4,\\
\label{eq:ida1b}
\vct{u}^{(n_{i})} =& \sum_{j=1}^{4} \mtx{T}^{(\nu_{3},ij)}\,\vct{v}^{(n_{j})},\qquad i = 1,\,2,\,3,\,4.
\end{eqnarray}
Observing that $\vct{u}^{(m_{2})} = \vct{u}^{(n_{n})}$ we combine (\ref{eq:ida1a}) for $i = 2$ with
(\ref{eq:ida1b}) for $i = 4$ to obtain the joint equation
\begin{multline}
\label{eq:ida2}
\mtx{T}^{(\nu_{1},21)}\,\vct{v}^{(m_{1})} +
\mtx{T}^{(\nu_{1},22)}\,\vct{v}^{(m_{2})} +
\mtx{T}^{(\nu_{1},23)}\,\vct{v}^{(m_{3})} +
\mtx{T}^{(\nu_{1},24)}\,\vct{v}^{(m_{4})}\\
=
\mtx{T}^{(\nu_{3},41)}\,\vct{v}^{(n_{1})} +
\mtx{T}^{(\nu_{3},42)}\,\vct{v}^{(n_{2})} +
\mtx{T}^{(\nu_{3},43)}\,\vct{v}^{(n_{3})} +
\mtx{T}^{(\nu_{3},44)}\,\vct{v}^{(n_{4})}.
\end{multline}
Utilizing further that $\vct{v}^{(m_{2})} = \vct{v}^{(n_{4})}$, we write
(\ref{eq:ida2}) along with (\ref{eq:ida1a}) and (\ref{eq:ida1b}) as
\begin{equation*}
\footnotesize
\left[\begin{array}{cccccc|c}
\mtx{T}^{(\nu_{1},11)} & \mtx{T}^{(\nu_{1},13)} & \mtx{T}^{(\nu_{1},14)} & 0 & 0 & 0 & \mtx{T}^{(\nu_{1},12)}\\
\mtx{T}^{(\nu_{1},31)} & \mtx{T}^{(\nu_{1},33)} & \mtx{T}^{(\nu_{1},34)} & 0 & 0 & 0 & \mtx{T}^{(\nu_{1},32)}\\
\mtx{T}^{(\nu_{1},41)} & \mtx{T}^{(\nu_{1},43)} & \mtx{T}^{(\nu_{1},44)} & 0 & 0 & 0 & \mtx{T}^{(\nu_{1},42)}\\
0 & 0 & 0 & \mtx{T}^{(\nu_{3},11)} & \mtx{T}^{(\nu_{3},12)} & \mtx{T}^{(\nu_{3},13)} & \mtx{T}^{(\nu_{3},14)}\\
0 & 0 & 0 & \mtx{T}^{(\nu_{3},21)} & \mtx{T}^{(\nu_{3},22)} & \mtx{T}^{(\nu_{3},23)} & \mtx{T}^{(\nu_{3},24)}\\
0 & 0 & 0 & \mtx{T}^{(\nu_{3},31)} & \mtx{T}^{(\nu_{3},32)} & \mtx{T}^{(\nu_{3},33)} & \mtx{T}^{(\nu_{3},34)}\\ \hline
 \mtx{T}^{(\nu_{1},21)} &  \mtx{T}^{(\nu_{1},23)} &  \mtx{T}^{(\nu_{1},24)} &
-\mtx{T}^{(\nu_{3},41)} & -\mtx{T}^{(\nu_{3},42)} & -\mtx{T}^{(\nu_{3},43)} &
 \mtx{T}^{(\nu_{1},22)} - \mtx{T}^{(\nu_{3},44)}
 \end{array}\right]\,
\left[\begin{array}{c}
\vct{v}^{(m_{1})} \\ \vct{v}^{(m_{3})} \\ \vct{v}^{(m_{4})} \\
\vct{v}^{(n_{1})} \\ \vct{v}^{(n_{2})} \\ \vct{v}^{(n_{3})} \\ \hline
\vct{v}^{(m_{2})}
\end{array}\right] =
\left[\begin{array}{c}
\vct{u}^{(m_{1})} \\ \vct{u}^{(m_{3})} \\ \vct{u}^{(m_{4})} \\
\vct{u}^{(n_{1})} \\ \vct{u}^{(n_{2})} \\ \vct{u}^{(n_{3})} \\ \hline
\vct{0}
\end{array}\right].
\end{equation*}
Eliminating $\vct{v}^{(m_{2})}$ from the system via a Schur complement yields the
operator $\mtx{T}^{(\sigma)}$ (upon suitable reblocking).

At this point, we have described a ``merge-two'' operation. A ``merge-four'' operation
can of course be obtained by simply combining three merge-two operations. To be precise,
suppose that $\tau$ is a node with the four children $\nu_{1},\,\nu_{2},\,\nu_{3},\,\nu_{4}$.
We introduce the two ``intermediate'' boxes $\sigma_{1}$ and $\sigma_{2}$ as shown in the
following figure:
\begin{center}
\setlength{\unitlength}{1mm}
\begin{picture}(65,15)
\put(00,00){\includegraphics[height=15mm]{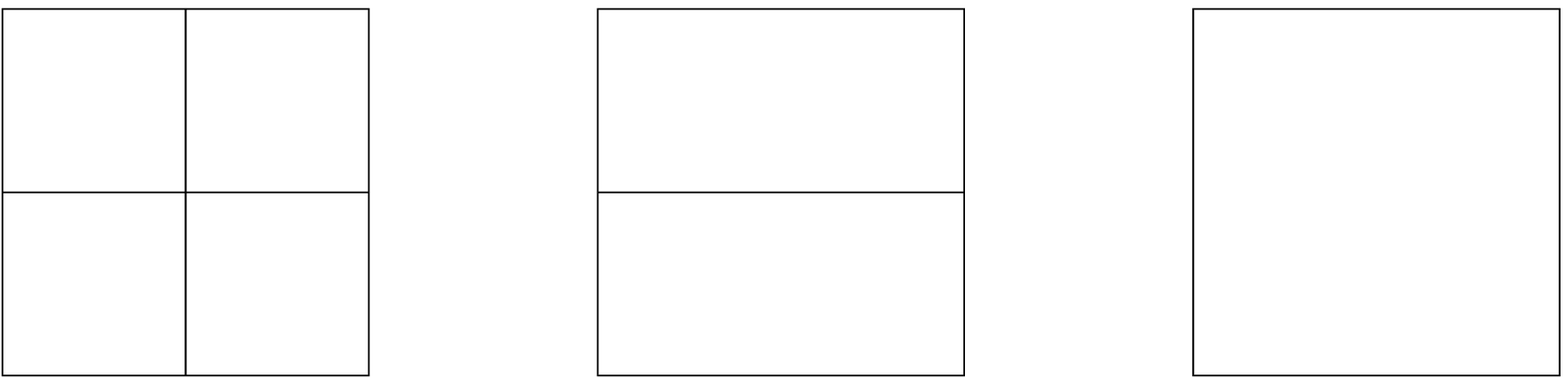}}
\put(02,03){$\nu_{1}$}
\put(02,10){$\nu_{2}$}
\put(09,03){$\nu_{3}$}
\put(09,10){$\nu_{4}$}
\put(18,07){$\Rightarrow$}
\put(30,03){$\sigma_{1}$}
\put(30,10){$\sigma_{2}$}
\put(41,07){$\Rightarrow$}
\put(54,07){$\tau$}
\end{picture}
\end{center}
Letting the procedure described earlier in the section be denoted by
``\texttt{merge$\underline{\mbox{ }}$two$\underline{\mbox{ }}$horizontal}''
and defining an analogous function
``\texttt{merge$\underline{\mbox{ }}$two$\underline{\mbox{ }}$vertical}''
we then find that
\begin{center}
\begin{tabular}{rl}
$\mtx{T}^{(\sigma_{1})} =$ & \texttt{merge$\underline{\mbox{ }}$two$\underline{\mbox{ }}$horizontal}$(\mtx{T}^{(\nu_{1})},\,\mtx{T}^{(\nu_{3})})$,\\
$\mtx{T}^{(\sigma_{2})} =$ & \texttt{merge$\underline{\mbox{ }}$two$\underline{\mbox{ }}$horizontal}$(\mtx{T}^{(\nu_{2})},\,\mtx{T}^{(\nu_{4})})$,\\
$\mtx{T}^{(\tau)} =$ & \texttt{merge$\underline{\mbox{ }}$two$\underline{\mbox{ }}$vertical}$(\mtx{T}^{(\sigma_{1})},\,\mtx{T}^{(\sigma_{2})})$.
\end{tabular}
\end{center}

\subsection{Fast construction of the N2D operator for a parent}
\label{sec:fast_merge}
The merge operation described in Section \ref{sec:slow_merge} has asymptotic cost $O(N^{1.5})$,
where $N$ is the total number of points on the edges of the leaves. To simplify slightly, the
reason is that forming the ``merge'' operation requires matrix inversion and matrix-matrix-multiplications
for of dense matrices whose size eventually grow to $O(\sqrt{N}) \times O(\sqrt{N})$.
However, these matrices all have internal structure. To be precise, in (\ref{eq:writtenout}),
the off-diagonal blocks have low rank, and the diagonal blocks are all Hierarchically Semi-Separable
(HSS) matrices. This means that accelerated matrix algebra can be used.
For a matrix of size $N' \times N'$, inversion can in fact be executed in $O(N')$ operations
(provided that the ``HSS-rank'' is a fixed low number, which it is in this case).

The acceleration procedure proposed here is analogous to the one described in \cite{2007_leborne_grasedyck,2010_ying_nested_diss,2007_gu_fem}.

\bibliography{main_bib}
\bibliographystyle{amsplain}

\clearpage

\begin{figure}
\setlength{\unitlength}{1mm}
\begin{picture}(50,50)
\put(000,000){\includegraphics[width=50mm]{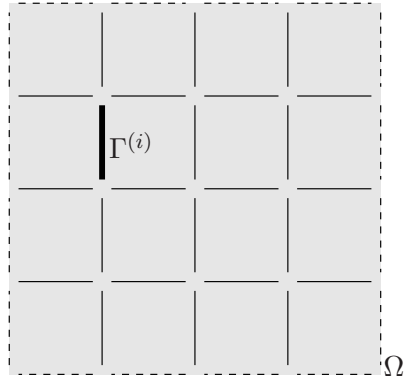}}
\put(013.5,029){$\Gamma^{(i)}$}
\put(50,00){$\Omega$}
\end{picture}
\caption{The computational box $\Omega$ (gray) is split into $16$ small boxes.
There are a total of $40$ edges in the discretization, 24 interior ones (solid lines)
and 16 exterior ones (dashed lines). One interior edge $\Gamma^{(i)}$ is marked with
a bold line. (Each edge continues all the way to the corner, but has been drawn slightly
forshortened for clarity.)}
\label{fig:domain}
\end{figure}

\begin{figure}
\setlength{\unitlength}{1mm}
\begin{picture}(50,50)
\put(000,000){\includegraphics[width=50mm]{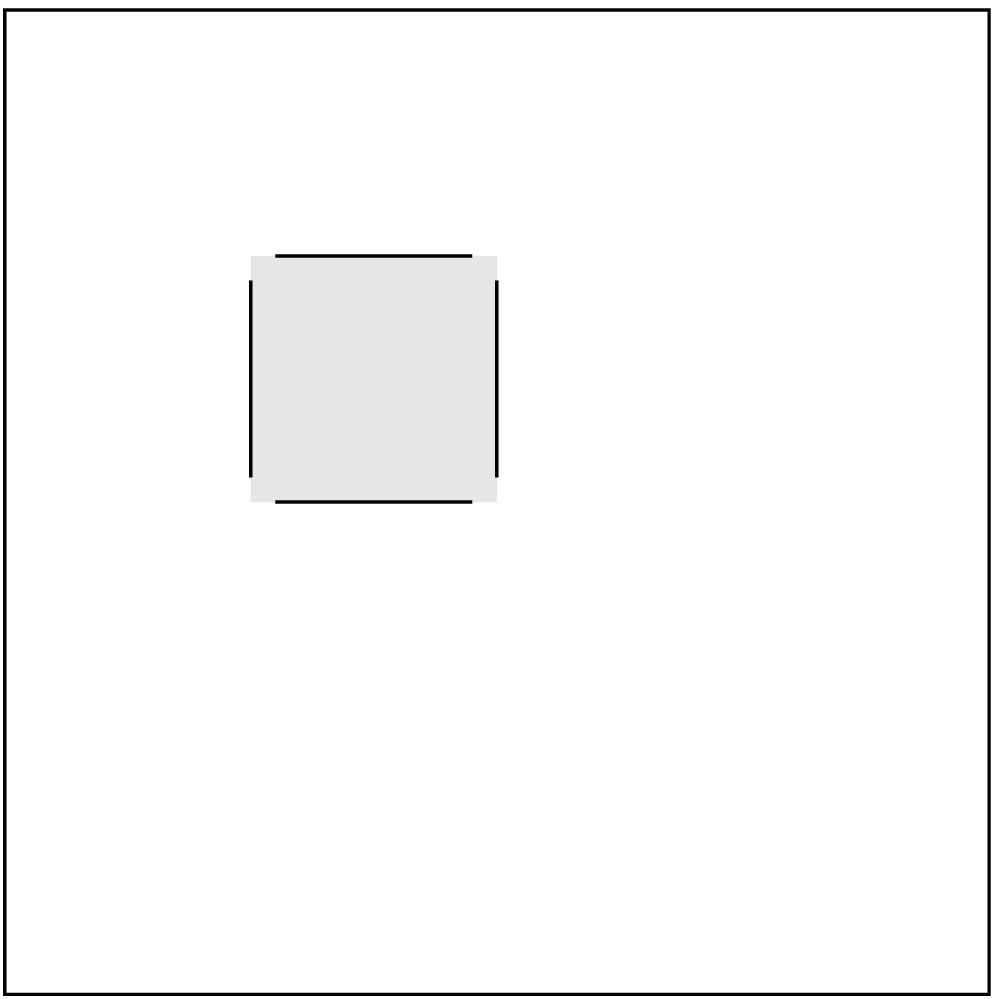}}
\put(016,030){$\Omega^{(\tau)}$}
\put(016,020){$\Gamma^{(i_{1})}$}
\put(026,030){$\Gamma^{(i_{2})}$}
\put(016,038){$\Gamma^{(i_{3})}$}
\put(005,030){$\Gamma^{(i_{4})}$}
\put(50,00){$\Omega$}
\end{picture}
\caption{The box $\Omega^{(\tau)}$ is marked in gray. Its edges are
$\Gamma^{(i_{1})},\,\Gamma^{(i_{2})},\,\Gamma^{(i_{3})},\,\Gamma^{(i_{4})}$.}
\label{fig:N2D}
\end{figure}

\begin{figure}
\setlength{\unitlength}{1mm}
\begin{picture}(50,50)
\put(000,000){\includegraphics[width=50mm]{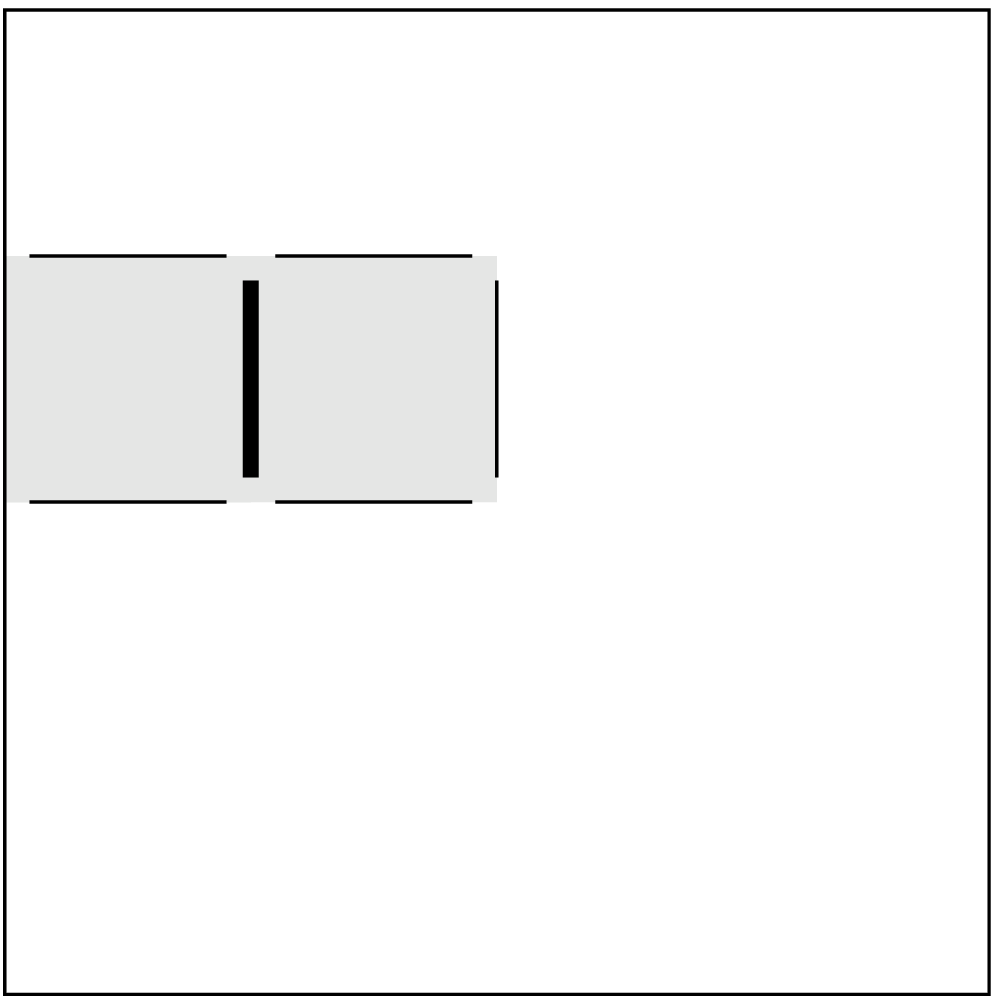}}
\put(016,030){$\Omega^{(\tau_{2})}$}
\put(003,030){$\Omega^{(\tau_{1})}$}
\put(003,020){$\Gamma^{(m_{1})}$}
\put(003,038){$\Gamma^{(m_{3})}$}
\put(-09,030){$\Gamma^{(m_{4})}$}
\put(016,020){$\Gamma^{(n_{1})}$}
\put(026,030){$\Gamma^{(n_{2})}$}
\put(016,038){$\Gamma^{(n_{3})}$}
\put(50,00){$\Omega$}
\end{picture}
\caption{Construction of the equilibrium equation for the edge $\Gamma^{(i)}$ in Figure \ref{fig:domain}.
It is the common edge of the boxes $\Omega^{(\tau_{1})}$ and $\Omega^{(\tau_{2})}$, which have
edges $\{\Gamma^{(m_{1})},\,\Gamma^{(m_{2})},\,\Gamma^{(m_{3})},\,\Gamma^{(m_{4})}\}$, and
      $\{\Gamma^{(n_{1})},\,\Gamma^{(n_{2})},\,\Gamma^{(n_{3})},\,\Gamma^{(n_{4})}\}$, respectively.
Observe that $\Gamma^{(i)} = \Gamma^{(m_{2})} = \Gamma^{(n_{4})}$ (the bold line).}
\label{fig:merge_new}
\end{figure}

\begin{figure}
\setlength{\unitlength}{1mm}
\begin{picture}(145,135)
\put(000,005){\includegraphics[width=145mm]{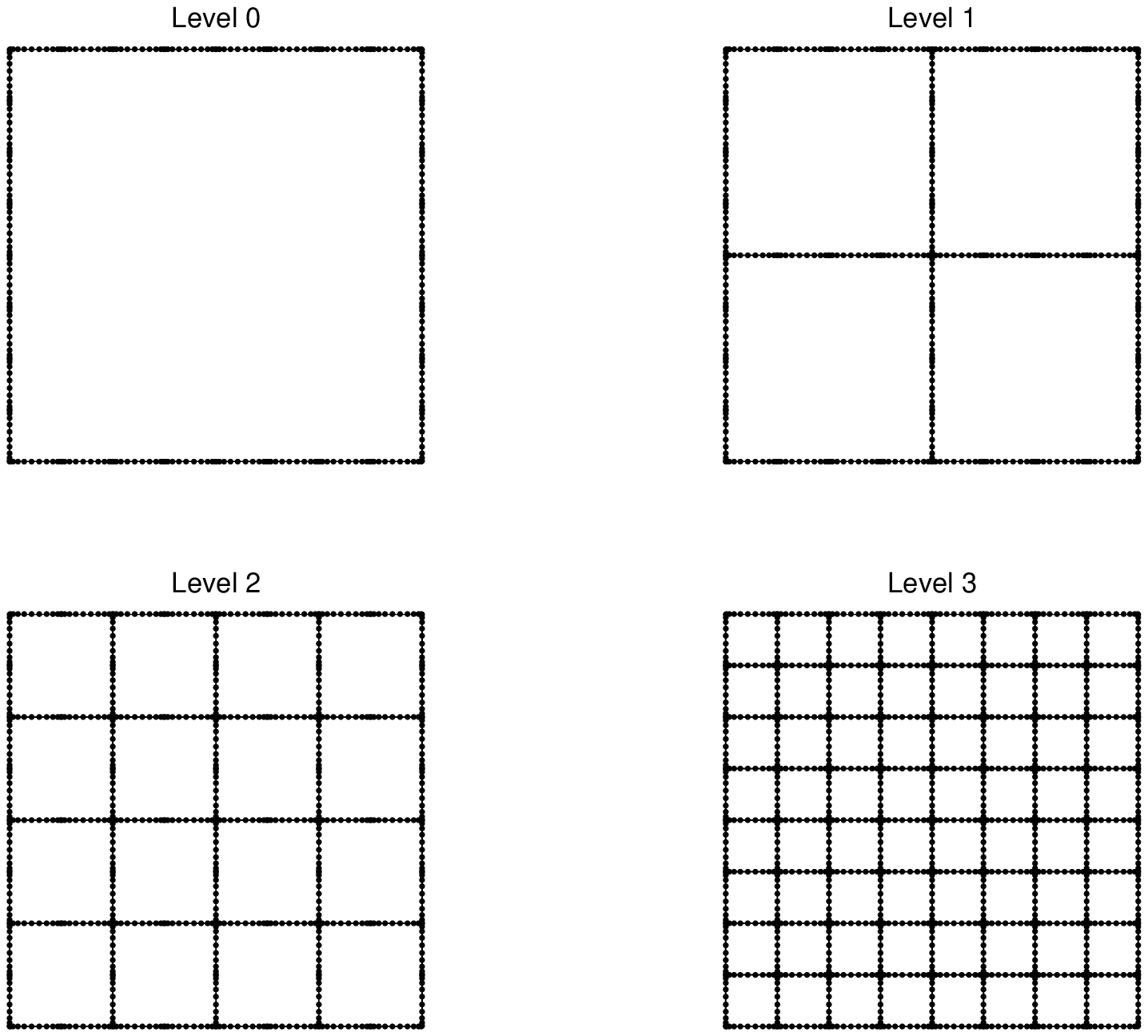}}
\put(020,115){$\tau=1$}
\put(095,085){$\tau=2$}
\put(095,115){$\tau=3$}
\put(125,085){$\tau=4$}
\put(125,115){$\tau=5$}
\put(004,010){$\tau=6$}
\put(004,025){$\tau=7$}
\put(016,010){$\tau=8$}
\put(016,025){$\tau=9$}
\put(023,071){(a)}
\put(114,071){(b)}
\put(023,000){(c)}
\put(114,000){(d)}
\end{picture}
\caption{Tree structure for a tree with $L = 3$ levels. There are
$10$ Gaussian nodes on each side of the leaf boxes. The black dots
mark the points at which the solution $\phi$ and its derivative (in the
direction normal to the indicated patch boundary) are tabulated.}
\label{fig:tree}
\end{figure}

\begin{figure}
\setlength{\unitlength}{1mm}
\begin{picture}(110,48)
\put(000,006){\includegraphics[width=110mm]{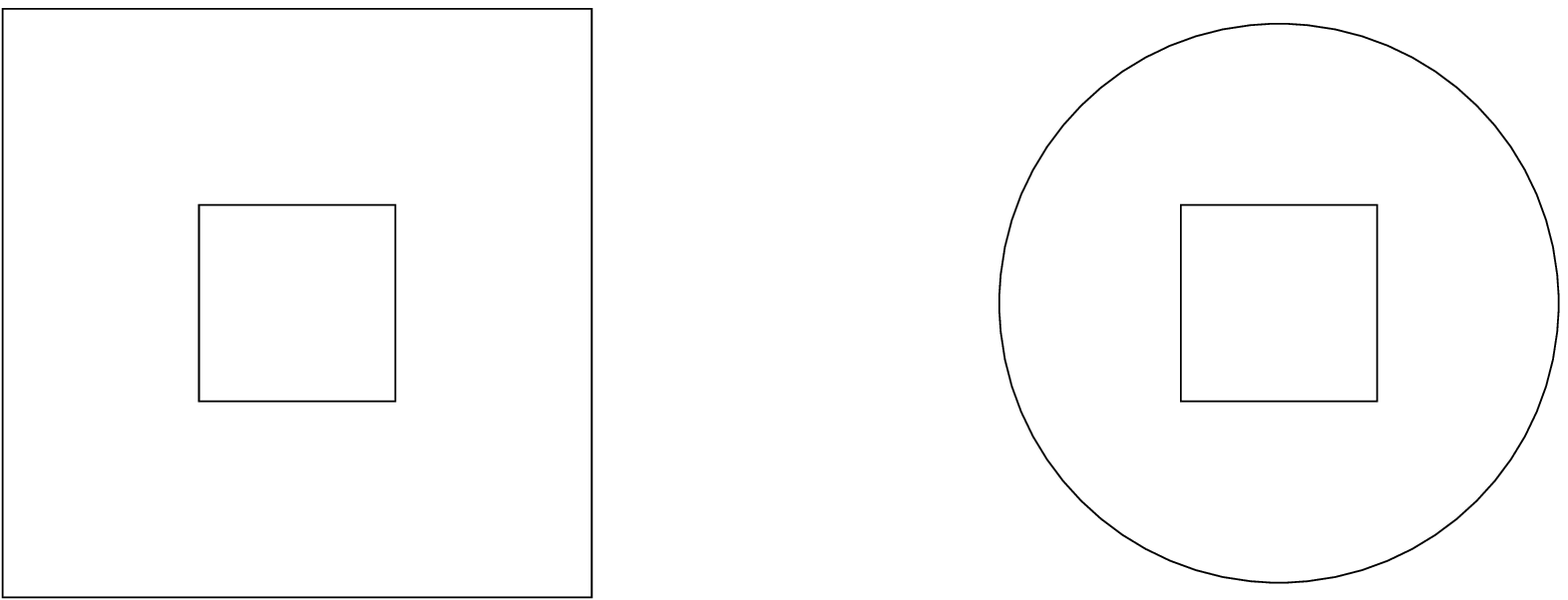}}
\put(018,000){(a)}
\put(090,000){(b)}
\put(018,029){$\Omega^{(\tau)}$}
\put(016,024){$\tilde{b} = b$}
\put(004,009){$\tilde{b}$ can be chosen freely}
\put(002,042){$\Psi$}
\put(088,029){$\Omega^{(\tau)}$}
\put(086,024){$\tilde{b} = b$}
\put(079,013){\footnotesize $\tilde{b}$ chosen freely}
\put(079,038){$\Psi$}
\end{picture}
\caption{Two choices of geometry for the local patch computation.
The choice (a) is natural since it conforms to the overall geometry.
The advantage of choice (b) is that the FFT can be used in the angular direction.}
\label{fig:local_patches}
\end{figure}

\begin{figure}
\setlength{\unitlength}{1mm}
\begin{picture}(110,65)
\put(000,000){\includegraphics[width=110mm]{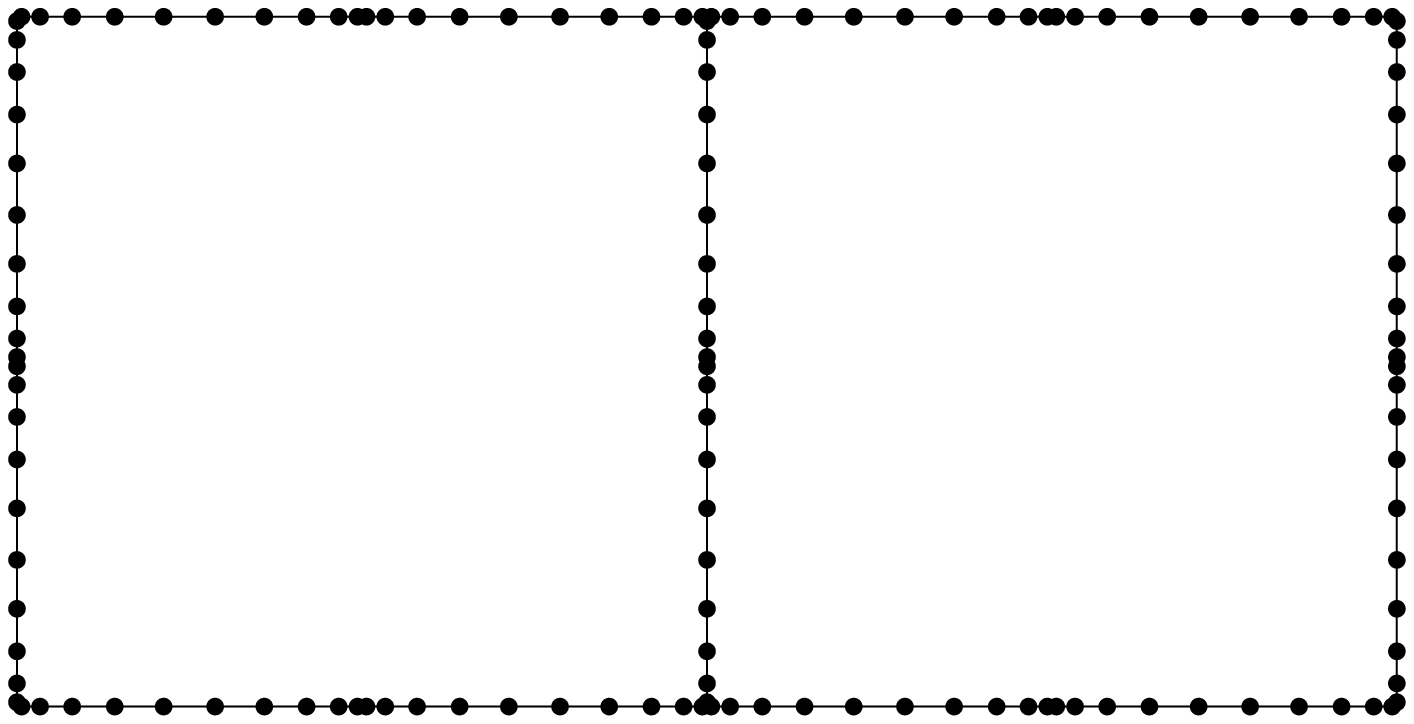}}
\put(026,033){$\Omega^{(\nu_{1})}$}
\put(077,033){$\Omega^{(\nu_{3})}$}
\put(005,028){$\Gamma^{(m_{4})}$}
\put(042,028){$\Gamma^{(m_{2})}$}
\put( 95,028){$\Gamma^{(n_{2})}$}
\put(058,028){$\Gamma^{(n_{4})}$}
\put(026,007){$\Gamma^{(m_{1})}$}
\put(026,050){$\Gamma^{(m_{3})}$}
\put(077,007){$\Gamma^{(n_{1})}$}
\put(077,050){$\Gamma^{(n_{3})}$}
\end{picture}
\caption{Geometry of the \textit{merge} operation.}
\label{fig:merge}
\end{figure}

\end{document}